\documentclass[leqn,9pt,a4paper]{amsart}
\usepackage{amsmath}
\usepackage[dvipdfmx]{graphicx}
\usepackage{multicolumn}
\usepackage{mikibase}
\usepackage{bm,color}
\usepackage[top=30truemm,bottom=30truemm,left=18truemm,right=18truemm]{geometry}

\begin{document}
\thispagestyle{empty}

\renewcommand{\thefootnote}{\fnsymbol{footnote}}
\begin{center}
{\bf \huge Numerical scheme based on the spectral method \\
for calculating nonlinear hyperbolic evolution equations}　
\end{center}　\vspace{10mm}　\\
\begin{multicols}{2}
\begin{center}
{\Large \bf Yoritaka Iwata}  \vspace{1.5mm}　\\
{ \bf  Kansai Uversity} \vspace{1.5mm} \\ 
{ \bf  Yamate-cho 3-3-35, Osaka 564-8680, Japan} \vspace{1.5mm} \\ 
{ \bf  (+81) 6-6368-1121} \vspace{1.5mm} \\ 
{ \bf  iwata$\_$phys@08.alumni.u-tokyo.ac.jp} \vspace{12mm} \\ 
\end{center}
\begin{center}
{\Large \bf Yasuhiro Takei}  \vspace{1.5mm} \\ 
{\bf Mizuho Information \& Research Institute Inc.}  \vspace{1.5mm}  \\
{ \bf  Kanda-Nishiki-cho 2-3, Tokyo 101-8443, Japan} \vspace{1.5mm} \\  
{ \bf  (+81) 3-5281-7500} \vspace{6mm} \\
\end{center} 
　\end{multicols}

\begin{multicols}{2}

\subsection*{\large ABSTRACT}
High-precision numerical scheme for nonlinear hyperbolic evolution equations is proposed based on the spectral method. The detail discretization processes are discussed in case of one-dimensional Klein-Gordon equations. In conclusion, a numerical scheme with the order of total calculation cost $O(N \log 2N)$ is proposed. As benchmark results, the relation between the numerical precision and the discretization unit size are demonstrated.  

\subsection*{\large CCS Concepts} 
 Mathematics of computing partial differential equations.

\subsection*{\large Keywords}
Fourier spectral method, high-precision calculation. \vspace{1.5mm} \\

\section{Introduction}
For concrete examples of hyperbolic evolution equations, one-dimensional linear and nonlinear Klein-Gordon equations are taken.
The initial and boundary values problem of one-dimensional Klein-Gordon equations 
\begin{equation} \begin{array}{ll} \label{eq01}
\ \frac{\partial^{2} u}{\partial t^{2}} + \alpha \frac{\partial^{2} u}{\partial x^{2}} + \beta F \left(u \right) = 0,  \vspace{3mm} \\
\ u(x, 0) = f(x), \quad
 \ \frac{\partial u}{\partial t}(x, 0) = g(x),  \vspace{3mm} \\
\ u(0, t) = u(L, t),  \quad
\ \frac{\partial u}{\partial t}(0, t) = \frac{\partial u}{\partial t}(L, t)
\end{array} \end{equation}
is considered for $(x,t) \in \left[0,L \right] \times \left[0, T \right]$, where $\alpha$, $\beta$ and $T$ are real numbers, and $f(x)$ and $g(x)$ are initial functions. 
The inhomogeneous term $F(u)$ is either linear or nonlinear function of $u$. 
This problem is also written by
\begin{equation} \begin{array}{ll} \label{eq02}
 \ \frac{\partial u}{\partial t} = v,  \vspace{3mm} \\ 
\ \frac{\partial v}{\partial t} + \alpha \frac{\partial ^{2} u}{\partial x^{2}} + \beta F \left(u \right) = 0,  \vspace{3mm} \\
\ u(x, 0) = f(x), \quad
 \ v(x, 0) = g(x),  \vspace{3mm} \\ 
\ u(0, t) = u(L, t),  \quad
\ v(0, t) = v(L, t).
\end{array} \end{equation}
In this manner, the first order evolution problem of hyperbolic type is obtained. 
This equation is regarded as wave equations.
Several conservative quantities in association with the wave propagation is utilized to confirm the precision of scheme.

In this article, for the initial and boundary values problem (\ref{eq02}), a precise numerical scheme is proposed. 
The Fourier spectral method is implemented for the spatial direction, and the $\theta$ scheme is introduced for the time direction. 
Note that $\theta =1/2$ is practically adopted in numerical benchmark tests, which is known as the Crank-Nicolson method.
Consequently the obtained numerical scheme with well-controlled precision is used for the finite-dimensional representation of infinite-dimensional dynamical systems.

\section{Discretization}
\subsection{Discretization of space}
The spectral method is employed [1]. 
The solution of (\ref{eq02}) is assumed to be expanded by the Fourier series:
\[ \begin{array}{ll} \label{eq03}
\ u(x, t) =  \\ \quad
a_{0}(t) + {\displaystyle \sum^{N}_{k=1 }} a_{k}(t) \cos \left( \frac{2 \pi}{L} k x \right) + {\displaystyle \sum^{N}_{k=1 }} b_{k}(t) \sin \left( \frac{2 \pi}{L} k x \right), \vspace{3mm} \\
 v(x, t) =\\ \quad
 c_{0}(t)  + {\displaystyle \sum^{N}_{k=1 }} c_{k}(t) \cos \left( \frac{2 \pi}{L} k x \right)  + {\displaystyle \sum^{N}_{k=1 }} d_{k}(t) \sin \left( \frac{2 \pi}{L} k x \right) .
\end{array} \]
Let us terminate the expansion by the $N$th term. 
Here, according to the terminology used in the mathematical physics, the space spanned by $x$ is referred to the coordinate space, and that spanned by $l$ to the momentum space. 
After substituting them to the 1st equation of (2), multiplying $\cos(2 \pi lx/L)$ and $\sin(2 \pi l x/L)$, and integrating by $x$ for $[0,L]$,
\begin{equation}
\label{eq:eq04}
\ \tfrac{d a_{0}}{d t} = c_{0}, \\
\ \tfrac{d a_{l}}{d t} = c_{l}, \\
\ \tfrac{d b_{l}}{d t} = d_{l}
\end{equation}
follow, where $l = 1, 2, \cdots , N$. 
In the same manner, after substituting them to the 2nd equation of (2), multiplying $\cos(2 \pi l x/L)$ and $\sin(2 \pi l x/L)$, and integrating by $x$ for $[0,L]$, we obtain
\begin{equation} \begin{array}{ll}
\label{eq:eq05}
L \tfrac{d c_{0}}{d t} + \beta \int^{L}_{0} F(u)  dx  = 0,  \vspace{3mm} \\
\frac{L}{2} \frac{d c_{l}}{d t} - \alpha \frac{2 \pi^{2}}{L}  l^{2}  a_{l} 
+ \beta \int^L_0 F(u) \cos ( \frac{2 \pi}{L}l x ) dx = 0, \vspace{3mm}\\
\frac{L}{2} \tfrac{d d_{l}}{d t} - \alpha \tfrac{2 \pi^{2}}{L}  l^{2}  b_{l} 
+ \beta \int^{L}_{0} F(u) \sin ( \tfrac{2 \pi}{L}l x )  dx = 0. \\
\end{array} \end{equation}
By solving Eqs. (3) and (4), the values of $a_0$, $c_0$, $a_l$, $b_l$, $c_l$, and $d_l$ are calculated. 
In terms of dealing with the nonlinearity, we pay attention to the integrals \\
\begin{itemize}
\item $\int_0^L F(u) dx$  \vspace{3mm} 
\item $\int_0^L F(u) \cos \left( \frac{2 \pi}{L} l x \right) dx$  \vspace{3mm} 
\item $\int_0^L F(u) \sin \left( \frac{2 \pi}{L} l x \right) dx$   \vspace{3mm}
\end{itemize}
in the right-hand side of (4). 
The operator-conversion method [1] is employed. 
At first, by introducing the Fourier inverse transformation, the nonlinear terms are separately calculated in the original coordinate space spanned by $x$. 
In the second, the nonlinear and linear terms are calculated in the momentum space spanned by $l$. 
This two-step treatment substantially reduces the calculation costs arising from the nonlinearity. 
In fact, using the Crank-Nicolson type numerical scheme, those integrals are approximated by
\begin{equation} \begin{array}{ll}
 \int^{L}_{0} F(u) \cos (\frac{2 \pi l x}{L }  ) dx 
 \simeq \frac{L}{J}  {\displaystyle \sum^{J-1}_{j=0}} \cos ( \frac{2 \pi l x_{j} }{L } ), \vspace{1mm} \\
 \int^{L}_{0} F(u) \sin ( \frac{2 \pi l x}{L}  ) dx 
 \simeq \frac{L}{J}  {\displaystyle \sum^{J-1}_{j=0}} \sin ( \frac{2 \pi l x_{j}}{L } ), \vspace{1mm} \\
 \int^{L}_{0} F(u) dx \simeq \frac{L}{J} {\displaystyle \sum^{J-1}_{j=0}} F(u_{j}). 
\end{array} \end{equation}
Under the periodic boundary condition, the spatial interval $[0,L]$ is equally discretized by $x_j$  with $j = 0, 1, \cdots , J$.
After the spatial discretization, unknown function $u$ at $t$ is denoted by $u_j = u(x_j, t)$. 
The similarity of the representation between the right-hand sides and the Fourier transform simplifies the calculations. Indeed, $u_j$ is calculated by the discrete Fourier transform (DFT), and $a_0$, $a_l$, and $b_l$ are calculated in the next.

In the present method, if $F(u)$ is $m$-th order polynomial, the left-hand sides and right-hand sides of  (5) coincide for $J \ge (M+1)N + 1$. 
By utilizing the Fast Fourier transform (FFT), the total calculation cost is reduced to $O(N \log 2 N)$. 
Consequently the original problem is spatially discretized as
 \begin{equation} \begin{array}{ll}
\label{eq:eq07}
\frac{d a_{0}}{d t} = c_{0}, \quad
\frac{d a_{l}}{d t} = c_{l}, \quad
\frac{d b_{l}}{d t} = d_{l},  \vspace{1.5mm} \\
L \frac{d c_{0}}{d t} + \beta \frac{L}{J} {\displaystyle \sum^{J-1}_{j=0}} F(u_{j})  = 0, \vspace{1.5mm} \\
\frac{L}{2} \frac{d c_{l}}{d t} - \frac{2 \pi^{2} l^2 \alpha  a_{l}}{L}   
+ \frac{L \beta}{J}  {\displaystyle \sum^{J-1}_{j=0}} F(u_{j}) \cos ( \frac{2 \pi l  x_{j}}{L}  ) = 0, \vspace{1.5mm} \\ \frac{L}{2} \frac{d d_{l}}{d t} -   \frac{2 \pi^{2} l^2 \alpha   b_{l}}{L} 
 + \frac{L \beta}{J}  {\displaystyle \sum^{J-1}_{j=0}} F(u_{j}) \sin (\frac{2 \pi l  x_{j}}{L}) = 0,
\end{array} \end{equation}
where $j = 0, 1, \cdots , J$ is the discretization in the coordinate space, and $l = 1, 2 \cdots , N$ is the discretization in the momentum space. 
Discretizations in two different spatial directions are mixed in this proposed formalism. 
In this sense, the spatial precision depends on both $J$ and $N$. 
Such a treatment is introduced for maintaining both the precision and the calculation cost.
Note that only momentum space discretization is enough for the linear cases.

\subsection{Discretization of time}
The $\theta$-scheme is employed [2]. 
For a sufficiently small $\Delta t$ and natural numbers $n$, we set $t_n = n \Delta t$. 
At first, using notations $a_l^n=a_l (t_n)$ and $c_l^n = c_l (t_n)$, the time discretization is given by
\[
\label{eq:eq08}
\begin{split}
\ \tfrac{a^{n+1}_{l} - a^{n}_{l}}{\Delta t} = \theta c^{n+1}_{l} &+ (1 - \theta) c^{n}_{l} ,
\end{split}
\]
where $\theta$ is a real number satisfying $0 \le \theta \le 1$. 
This is equivalent to
\[
\label{eq:eq09}
\begin{split}
\ a^{n+1}_{l} &= a^{n}_{l} + \Delta t [ (1-\theta) c^{n}_{l}
 +  \theta c^{n+1}_{l} ], 
\end{split}
\]
and similarly
\[
\label{eq:eq10}
\begin{split}
b^{n+1}_{l} &= b^{n}_{l} + \Delta t [ (1-\theta) d^{n}_{l}
 + \theta d^{n+1}_{l} ] 
\end{split}
\]
with a notation $b_l^n=b_l (t_n)$ and $d_l^n=d_l (t_n)$. 
In the second, using notations $u_j^n=u_l (x_j,t_n)$, and ${\hat F}_l^n=u_l(x_j,t_n )$, the fifth equation of (6) become
\begin{equation}
\label{eq:eq11}
\begin{split}
& \tfrac{L}{2} \tfrac{ (c^{n+1}_{l} - c^{n}_{l} )}{\Delta t} - \tfrac{2 \alpha \pi^{2}  l^{2}}{L} \theta a^{n+1}_{l} - \tfrac{2 \alpha \pi^{2}  l^{2} }{L}(1 - \theta) a^{n}_{l} \\
& \ \  + \beta \tfrac{L}{2} \theta \hat{F}^{n+1}_{l} 
+ \beta \tfrac{L}{2} (1 - \theta) \hat{F}^{n}_{l} = 0 
\end{split}
\end{equation}
which is equivalent to
\begin{equation} \begin{array}{ll}
\label{eq:eq12}
c^{n+1}_{l} = c^{n}_{l} + \alpha ( \tfrac{2 \pi l}{L} )^2 \Delta t [ (1 - \theta)  a^{n}_{l}
   +  \theta a^{n+1}_{l} ]  \vspace{1.5mm} \\
  - \beta \Delta t [ (1 - \theta)  \hat{F}^{n}_{l}
  + \theta   \hat{F}^{n+1}_{l} ], 
\end{array} \end{equation}
and in the same way
\begin{equation}
\label{eq:eq13}  \begin{array}{ll}
d^{n+1}_{l} = d^{n}_{l} 
+ \alpha ( \tfrac{2 \pi l}{L} )^2  \Delta t [ (1 - \theta) b^{n}_{l}
 + \theta  b^{n+1}_{l} ] \vspace{1.5mm} \\
  - \beta \Delta t [ (1 - \theta)   \hat{G}^{n}_{l}
+  \theta  \hat{G}^{n+1}_{l}]. \ 
\end{array} \end{equation}
Consequently the original problem is fully discretized as
\begin{equation} \begin{array}{ll}
\label{eq:eq15}
a^{n+1}_{l} = a^{n}_{l} + \Delta t (1-\theta) c^{n}_{l} + \Delta t \theta c^{n+1}_{l},  \vspace{3mm}  \\ 

b^{n+1}_{l} = b^{n}_{l} + \Delta t (1-\theta) d^{n}_{l} + \Delta t \theta d^{n+1}_{l},  \vspace{3mm}  \\ 

c^{n+1}_{l} = c^{n}_{l} + (1 - \theta)  \Delta t [ \alpha (\tfrac{2 \pi l}{L} )^2  a^{n}_{l}
- \beta  \hat{F}^{n}_{l} ] \\
\hspace{10mm} +  \theta \Delta t [ \alpha ( \tfrac{2 \pi l}{L} )^2 a^{n+1}_{l}
- \beta \hat{F}^{n+1}_{l} ],  \vspace{3mm} \\ 

d^{n+1}_{l} = d^{n}_{l} +  (1 - \theta) \Delta t [ \alpha (\tfrac{2 \pi l}{L})^2 b^{n}_{l}
- \beta  \hat{G}^{n}_{l} ] \\
\hspace{10mm} +  \theta \Delta t [ \alpha ( \tfrac{2 \pi l}{L} )^2 b^{n+1}_{l}
- \beta  \hat{G}^{n+1}_{l} ],  \vspace{3mm}  \\ 

c^{n+1}_{0} = c^{n}_{0}  - \beta (1 - \theta) \Delta t \hat{F}^{n}_{0}
 - \beta \theta \Delta t \hat{F}^{n+1}_{0} \\
\end{array} \end{equation}
in terms of time and space.  Solving (\ref{eq:eq15}) means calculating unknown functions $a^{n+1}_{0}$, $c^{n+1}_{0}$, $a^{n+1}_{l}$, $ b^{n+1}_{l}$, $c^{n+1}_{l}$, $d^{n+1}_{l}$ from the known functions $a^{n}_{0}$, $c^{n}_{0}$, $a^{n}_{l}$, $b^{n}_{l}$, $c^{n}_{l}$, $d^{n}_{l}$.
In general, the unknown functions are existing in the right-hand side of (\ref{eq:eq13}), so that it is necessary to be solved by the iterative scheme. 

\begin{figure*}[tb]
\begin{center}
  \includegraphics[width=7cm,clip]{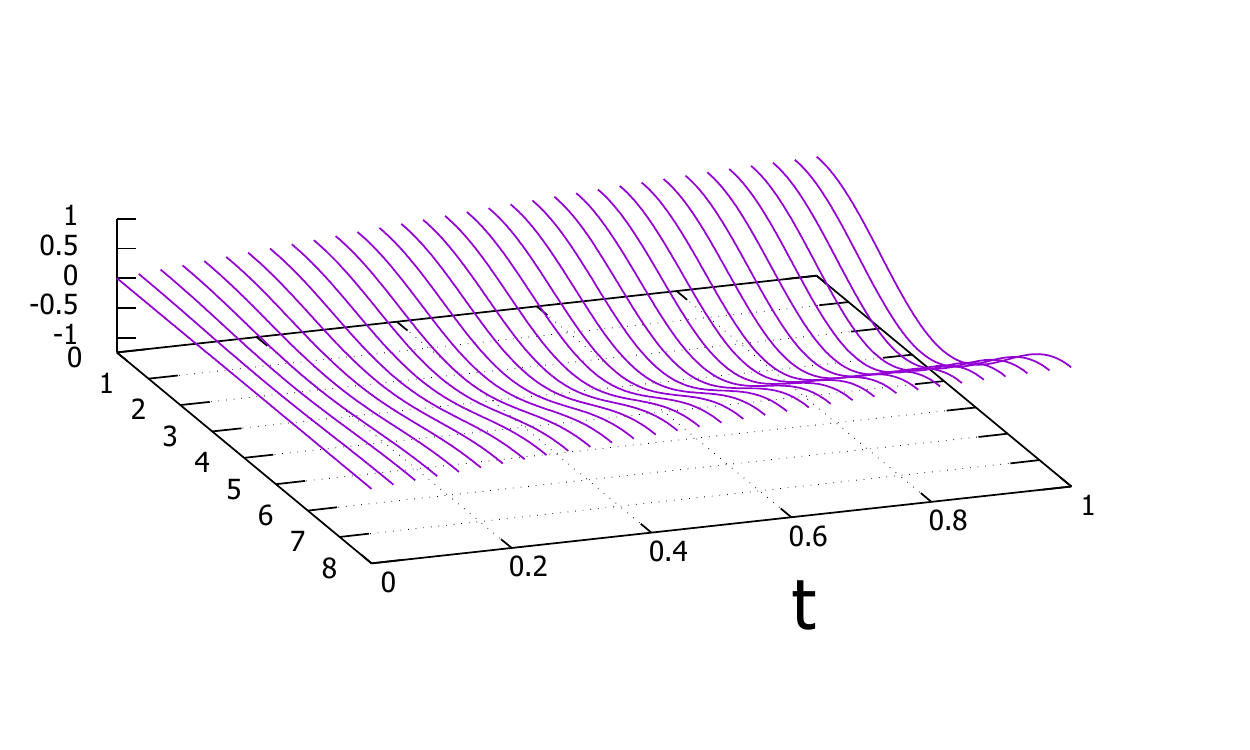}
  \includegraphics[width=7cm,clip]{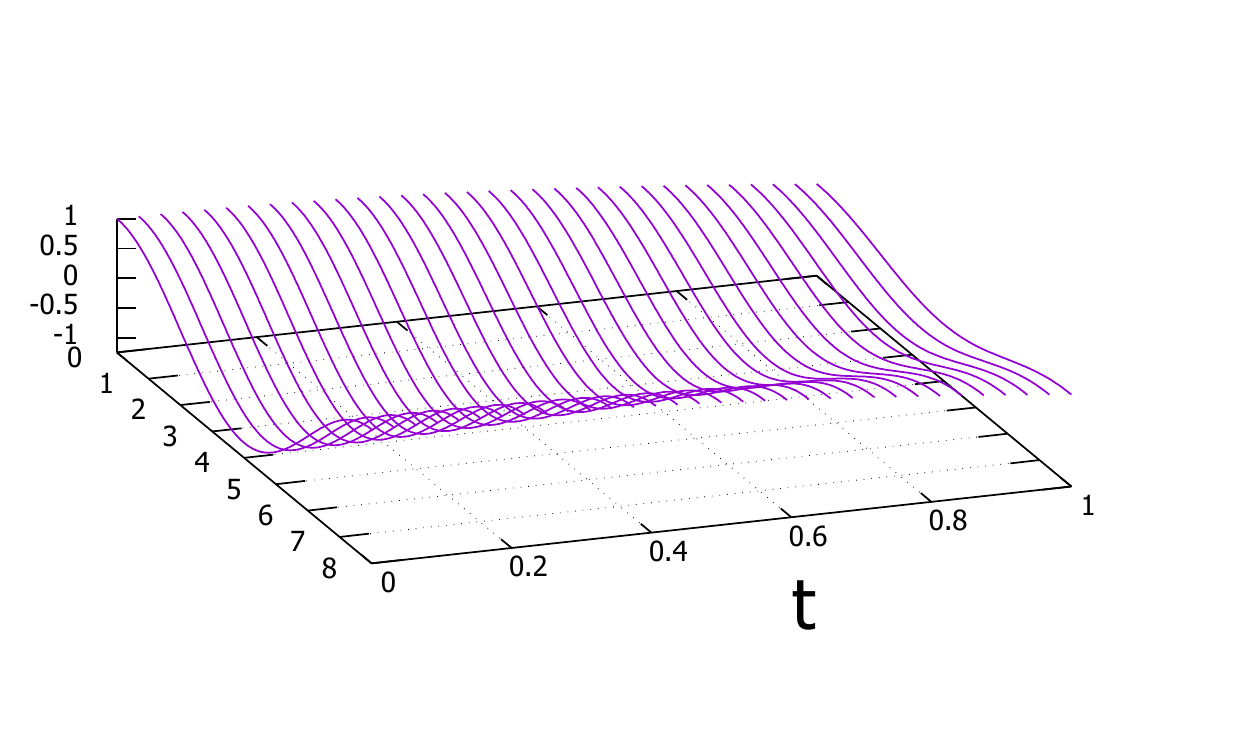}
\end{center}
\vspace{-7mm}
\begin{minipage}{14.5cm}
 \fgcaption{Linear Klein-Gordon dynamics: time evolution of $u$ and $v$.}
 \label{pic:pic02a}
 \vspace{1em}
\end{minipage}
\begin{center}
  \includegraphics[width=7cm,clip]{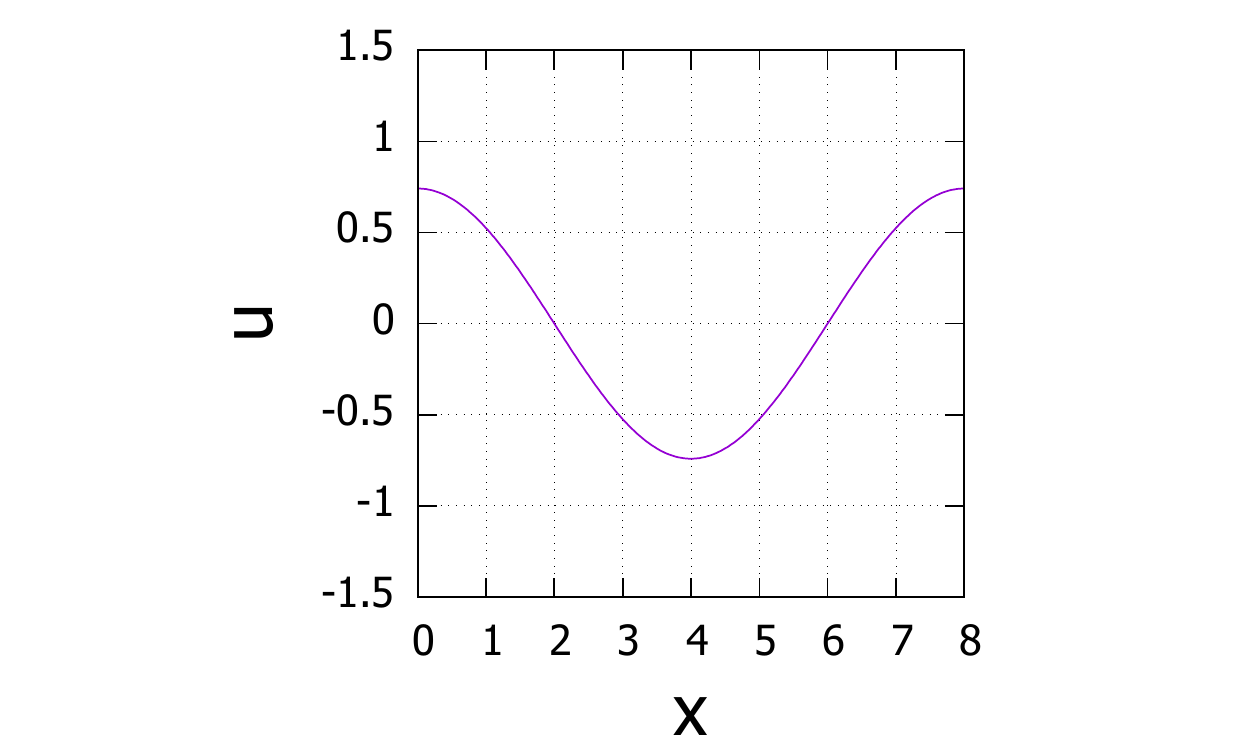}
  \includegraphics[width=7cm,clip]{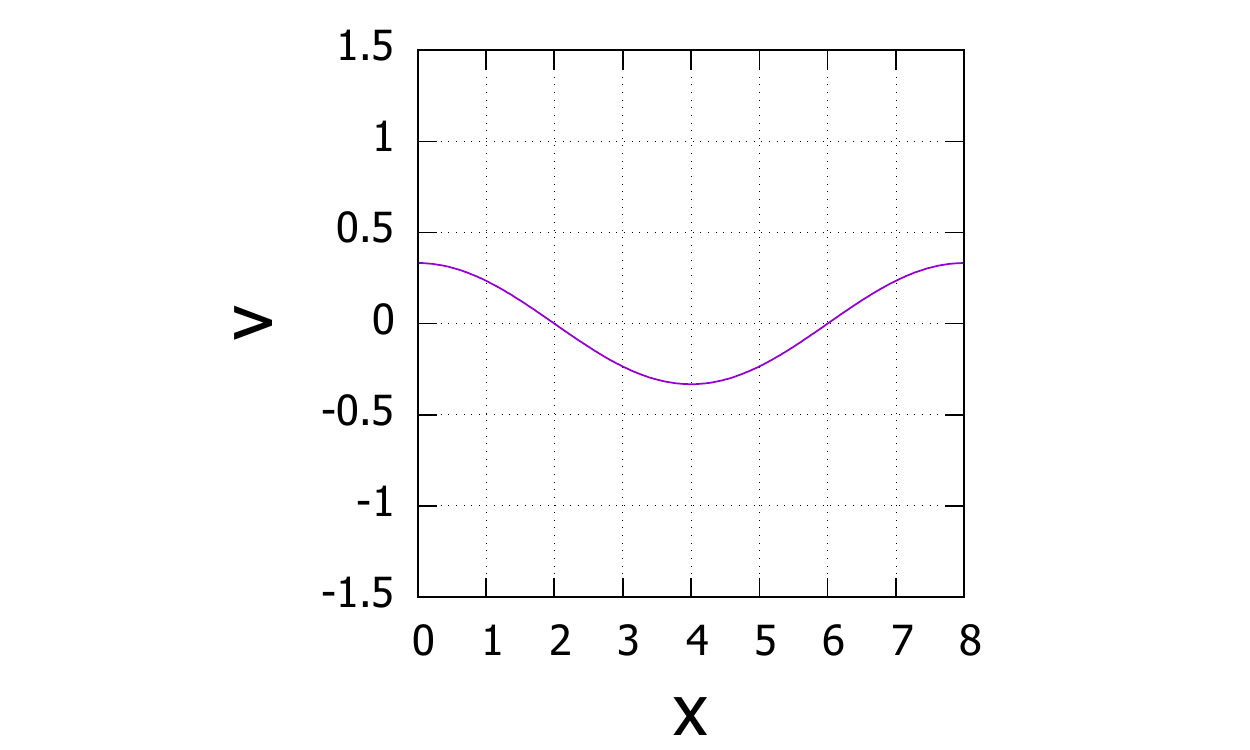}
\end{center}
\vspace{-2mm} 
\begin{minipage}{14cm}
 \fgcaption{Linear Klein-Gordon equation: $u$ and $v$ at $t = 1$.}
 \label{pic:pic04a}
\end{minipage}
\end{figure*}

\section{Iterative scheme}
The $\theta$ scheme solver is generally an implicit method. That is, the iterative treatment is required. 
Here we introduce additional parameter $\nu$ for the implicit treatment.  
Using a new parameter $\nu$, the intermediate stage is introduced in which $a^{n+1}_{0}$, $c^{n+1}_{0}$, $a^{n+1}_{l}$, $b^{n+1}_{l}$, $c^{n+1}_{l}$, $d^{n+1}_{l}$
are separated into 
$a^{n+1, \nu+1}_{0}$, $c^{n+1, \nu+1}_{0}$, $a^{n+1, \nu+1}_{l}$, $b^{n+1, \nu+1}_{l}$, $c^{n+1, \nu+1}_{l}$, $d^{n+1, \nu+1}_{l}$ and
$a^{n+1, \nu}_{0}$, $c^{n+1, \nu}_{0}$, $a^{n+1, \nu}_{l}$, $b^{n+1, \nu}_{l}$, $c^{n+1, \nu}_{l}$, $d^{n+1, \nu}_{l}$, and it follows that 
\begin{equation}
\label{eq:eq16} \begin{array}{ll}
a^{n+1, \nu+1}_{l} = a^{n}_{l} + \Delta t (1-\theta) c^{n}_{l}
+ \Delta t \theta c^{n+1, \nu}_{l} ,  \vspace{3mm} \\  

b^{n+1, \nu+1}_{l} = b^{n}_{l} + \Delta t (1-\theta) d^{n}_{l}
+ \Delta t \theta d^{n+1, \nu}_{l}, \vspace{3mm} \\

c^{n+1, \nu+1}_{l} = c^{n}_{l}
 +  (1 - \theta) \Delta t [ \alpha ( \tfrac{2 \pi l}{L} )^2 a^{n}_{l}
- \beta \hat{F}^{n}_{l} ] \vspace{1.5mm} \\
\hspace{10mm} + \theta \Delta t [ \alpha (\tfrac{2 \pi l}{L} )^2 a^{n+1, \nu}_{l}
- \beta \hat{F}^{n+1, \nu}_{l} ],  \vspace{3mm}  \\ 

d^{n+1, \nu+1}_{l} = d^{n}_{l}
 +  (1 - \theta) \Delta t [ \alpha ( \tfrac{2 \pi l}{L} )^2 b^{n}_{l}
 - \beta  \hat{G}^{n}_{l} ] \vspace{1.5mm} \\
\hspace{10mm} +  \theta \Delta t [ \alpha ( \tfrac{2 \pi l}{L} )^2 b^{n+1, \nu}_{l}
 - \beta  \hat{G}^{n+1, \nu}_{l} ], \vspace{3mm} \\ 

c^{n+1, \nu+1}_{0} = c^{n}_{0}  - \beta (1 - \theta) \Delta t \hat{F}^{n}_{0}
 - \beta \theta \Delta t \hat{F}^{n+1, \nu}_{0}.
\end{array}
\end{equation}
Using the recurrence relation (\ref{eq:eq16}), the $(\nu+1)$th values $a^{n+1, \nu+1}_{0}$, $c^{n+1, \nu+1}_{0}$, $a^{n+1, \nu+1}_{l}$, $b^{n+1, \nu+1}_{l}$, $c^{n+1, \nu+1}_{l}$, $d^{n+1, \nu+1}_{l}$ are obtained by the $\nu$th values $a^{n+1, \nu}_{0}$, $c^{n+1, \nu}_{0}$, $a^{n+1, \nu}_{l}$, $b^{n+1, \nu}_{l}$, $c^{n+1, \nu}_{l}$, $d^{n+1, \nu}_{l}$.
This process is sequentially repeated until the convergence. 
The converged ones satisfy (10) in a numerical sense, then they are renamed as  $a^{n+1}_{0}$, $c^{n+1}_{0}$, $a^{n+1}_{l}$, $b^{n+1}_{l}$, $c^{n+1}_{l}$, $d^{n+1}_{l}$.
Here $\hat{F}^{n+1, \nu}_{l}$ and  $\hat{G}^{n+1, \nu}_{l}$ are obtained by
\begin{equation}
\label{eq:eq17}
\begin{array}{ll}
\hat{F}^{n+1, \nu}_{0} = \tfrac{1}{J} {\displaystyle \sum^{J-1}_{j=0}} F(u^{n+1, \nu}_{j}), \vspace{0.75em}\\
\hat{F}^{n+1, \nu}_{l} = \tfrac{2}{J}  {\displaystyle  \sum^{J-1}_{j=0}} F(u^{n+1, \nu}_{j}) 
 \cos ( \tfrac{2 \pi}{L} l x_{j} ), \vspace{0.75em}\\
\hat{G}^{n+1, \nu}_{l} = \tfrac{2}{J}  {\displaystyle  \sum^{J-1}_{j=0}} F(u^{n+1, \nu}_{j}) 
 \sin ( \tfrac{2 \pi}{L} l x_{j} ).
\end{array}
\end{equation}
In particular $u_j^{n+1,\nu}$ being necessary at the intermediate stage is obtained by the DFT in which $a^{n+1, \nu}_{l}$, $b^{n+1, \nu}_{l}$, $c^{n+1, \nu}_{l}$, and $d^{n+1, \nu}_{l}$ are utilized.

The iteration process is the bottle neck of the calculation. 
Consequently the total calculation cost is mostly dominated by the integration of (\ref{eq:eq15}) in the iteration process, whose order is $O(N \log_{2}N)$. 
The order of cost of our total calculation is concluded to be  $O(N \log_{2}N)$.

\begin{figure*}[tb]
\vspace{-10mm} 
\begin{center}
  \includegraphics[width=7cm ,clip]{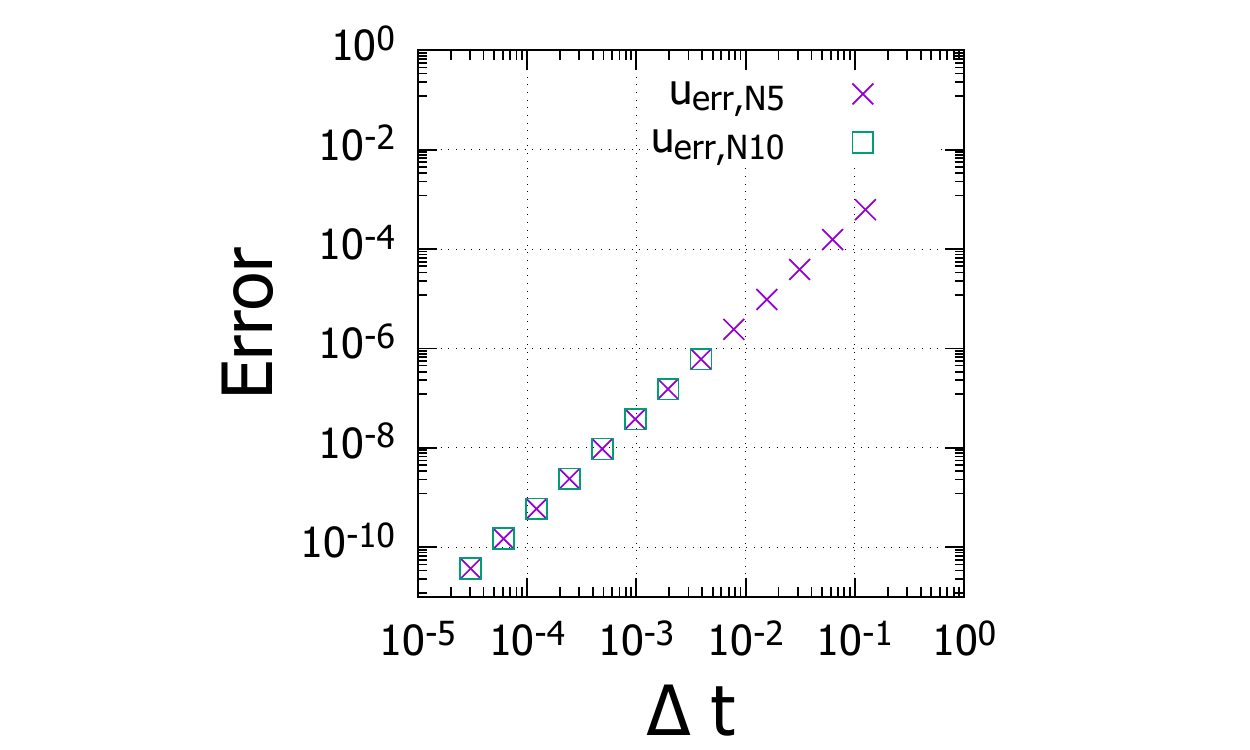}
  \includegraphics[width=7cm ,clip]{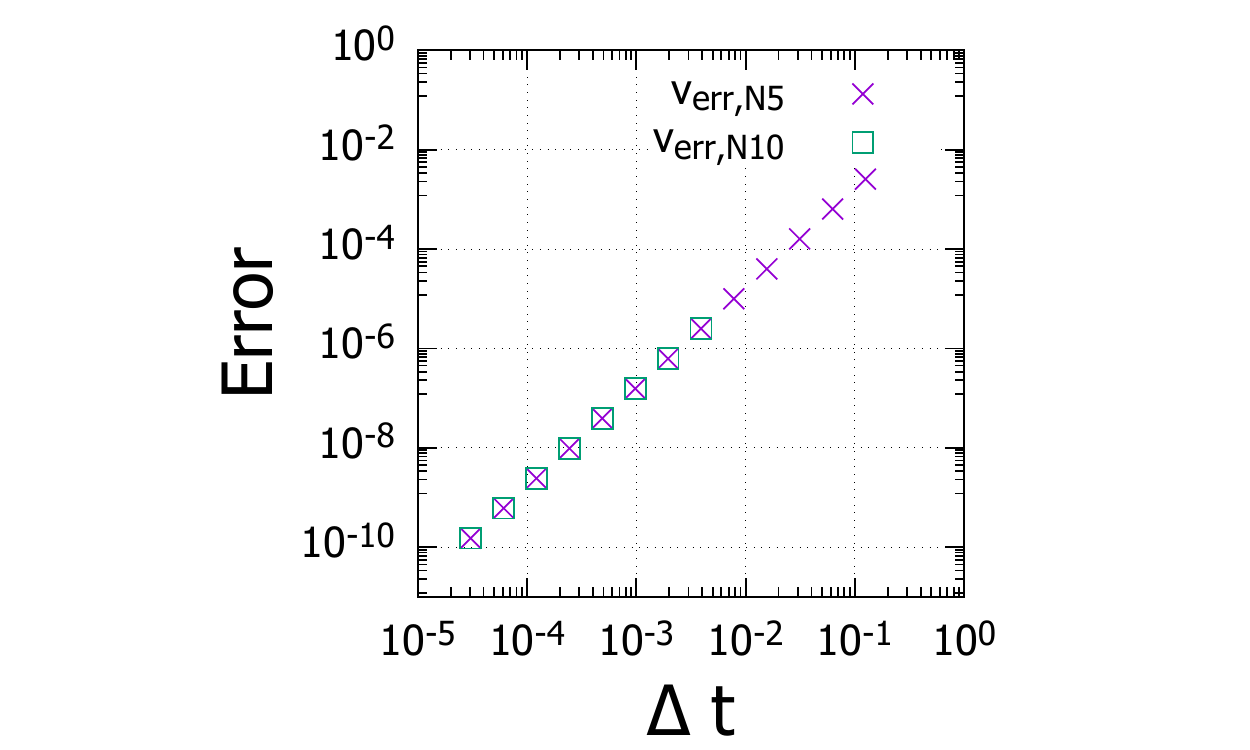}
\vspace{-2mm} 
\fgcaption{Linear case: the error dependence on $\Delta t$. Two different cases $N=2^5$ (labeled by $u_{err}\_N5$ and $v_{err}\_N5$) and $N=2^{10}$  (labeled by $u_{err}\_N10$ and $v_{err}\_N10$) are examined.}
\label{pic:pic06}
\end{center}
\begin{center}
  \includegraphics[width=7cm ,clip]{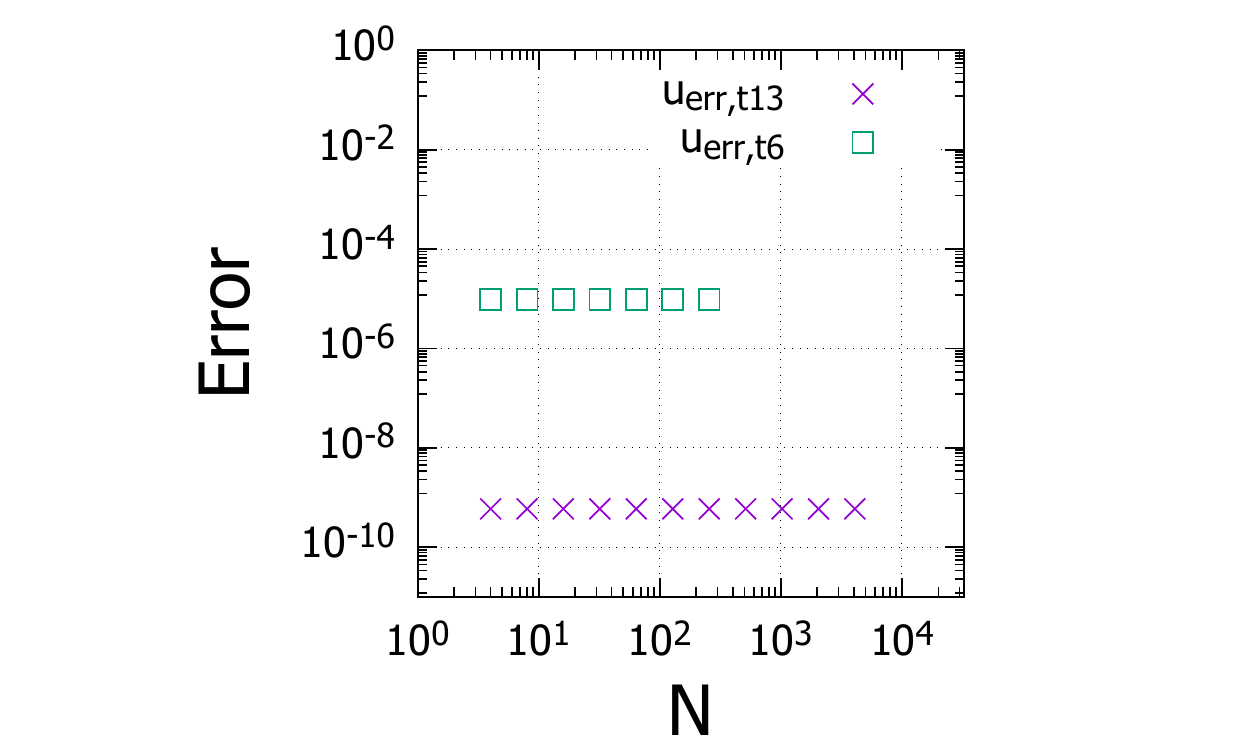}
  \includegraphics[width=7cm ,clip]{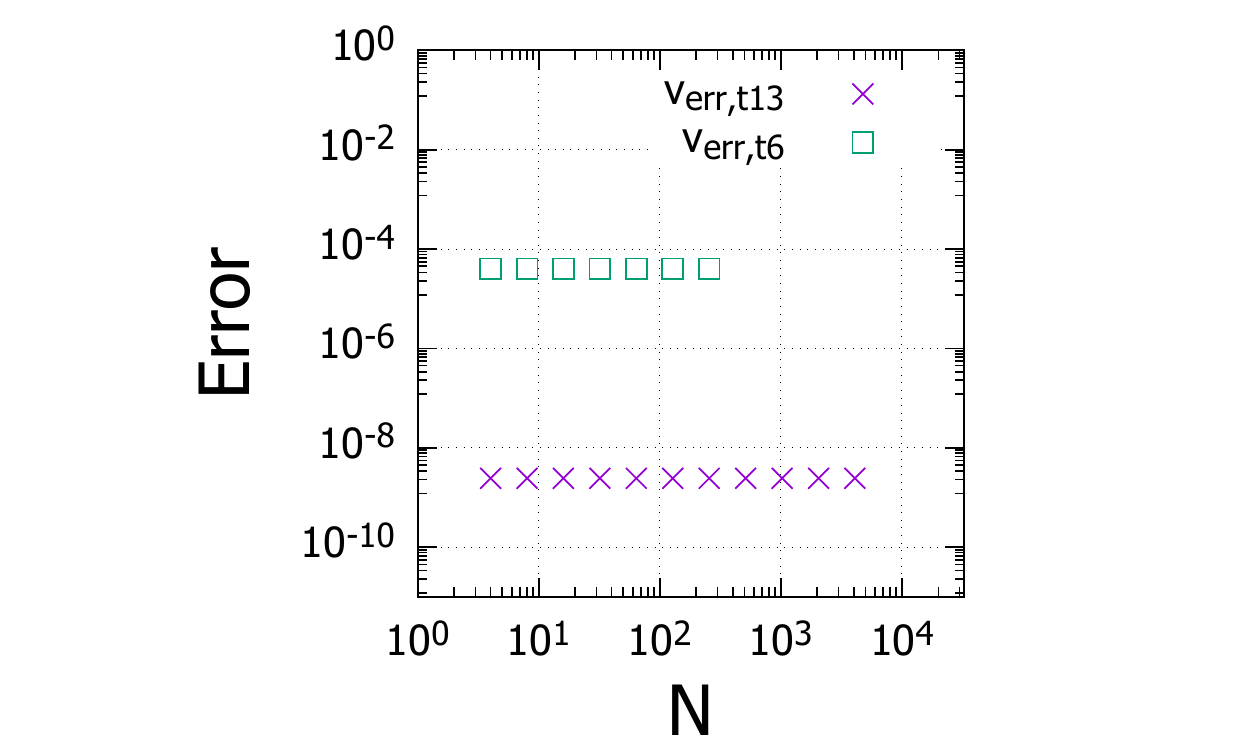}
\vspace{-2mm} 
\fgcaption{Linear case: the error dependence on $\Delta t$. Two different cases $t=2^{-6}$ (labeled by $u_{err}\_t6$ and $v_{err}\_t6$) and $t=2^{-13}$  (labeled by $u_{err}\_t13$ and $v_{err}\_t13$) are examined. }
\label{pic:pic08}
\end{center}
\end{figure*}

\section{Benchmark calculations with numerical error tests}
\subsection{Linear case}
For the initial and boundary values problem (\ref{eq02}),
\[ \begin{array}{ll}
F(u)=u, \quad \alpha = -1, \quad
\beta = 1,    \quad
\Omega = [0, L].
\end{array} \]
The corresponding problem is written by
\begin{equation}
\label{eq:eq18}
\begin{array}{ll}
\ \tfrac{\partial v}{\partial t} -  \tfrac{\partial ^{2} u}{\partial x^{2}} + u = 0, \vspace{3mm} \\
\ \tfrac{\partial u}{\partial t} = v, \vspace{3mm} \\ 
\ u(x, 0) = 0,  \quad  u(0, t) = u(L, t),\vspace{3mm} \\
\ v(x, 0) = \cos ( \tfrac{2 \pi}{L} x ),  \quad v(0, t) = v(L, t). \\ 
\end{array}
\end{equation}
An exact solution of this problem is given by
\begin{equation}
\label{eq:eq19}
\begin{array}{ll}
 u(x, t) =
 \frac{1}{\omega} \sin ( \omega t )
 \cos \left( \tfrac{2 \pi}{L} x \right), \vspace{0.75em}\\ 
 v(x, t) = \cos (\omega t ) 
 \cos \left( \tfrac{2 \pi}{L} x \right),  \vspace{0.75em} \\
\omega =  {\scriptstyle \sqrt{1 + \left( \tfrac{2 \pi}{L} \right)^{2}}},
\end{array}
\end{equation}
and the corresponding numerical solution is shown in Figs.~1 and 2.
In the actual numerical calculations, we further assume
\[ \begin{array}{ll}
\theta = \frac{1}{2},  \vspace{0.75em} \\ 
J \geq 2N + 1,  \vspace{0.75em} \\
L=8.
\end{array}  \]
The error comparing the exact and numerical solutions are shown in Figs.~ 3 and 4, where $N = 2^5$ and $2^{10}$ cases with $\Delta t = 2^{-2}$, $2^{-3}$, $\cdots$, $2^{15}$ cases are examined. 
In order to obtain the error estimates, the difference between the exact and numerical solutions is calculated at each point $x_j$.

After calculating absolute and relative errors at each point, the minimum error between absolute and relative errors are obtained at each point $x_j$. 
By denoting this minimum as ${\rm err}(x_j)$, our error function used in the plots is defined by
\[
Error = \max_{x_j} {\rm err}(x_j),
\]
where, in terms of picking up the error without having a overflow of number representation, the minimum is taken before taking the maximum. 
Two kinds of errors are used for calculating ${\rm err}(x_j),$ the relative error and the absolute error.
The relative error is usually adopted, and the absolute error is adopted  only when the absolute value of $u_i$ is too small for the value of relative error to be handled in the computer. 
Finally the maximum of error is searched over all the discretized points.

For the error arising from the time discretization, if we apply half of $\Delta t$, it results in the quartered error. 
It simply corresponds to the fact that the present scheme is 2nd order scheme for time direction. 
That is, it is possible to reduce/control the error rather easily by taking sufficiently small $\Delta t$.
This information is practical to determine the value of $\Delta t$.

Let us compare $N$ = $2^5$ case to $2^{10}$ case. 
In this comparison, there is no significant difference in error values at least if the same value is chosen for $\Delta t$. 
Note here that the convergence was turned out to be false when $N = 2^{10}$ with larger $\Delta t$ was applied. 
This remarkable property is an advantage of spectral method in which the error does
not depend on $N$ but on the $\Delta t$. 
This fact is also true in the other plots (Fig.~4). 
In Fig.~4, the error included in the numerical solution at $t = 1$ is shown, and the $N$-independence of error is clearly seen.
Indeed, in the present numerical scheme, remarkable properties
\begin{itemize}
\item if $F(u)$ is $M$th-order polynomial of $u$, the left- and right-hand sides of Eq. (4) being used in the calculation of $F(u)$ are exactly the same if $J \ge (M + 1)N + 1$ is satisfied; \\
\item for a sufficiently large $N$, the linear solution with a chosen initial value is exactly represented by the Fourier expansion;  \vspace{0.5mm}
\end{itemize}
are confirmed.
These contribute to the $N$ independence of possible errors.
In conclusion, the error of the present scheme can be controlled by the time discretization.

\begin{figure*}[tb]
\begin{center}
  \includegraphics[width=7cm,clip]{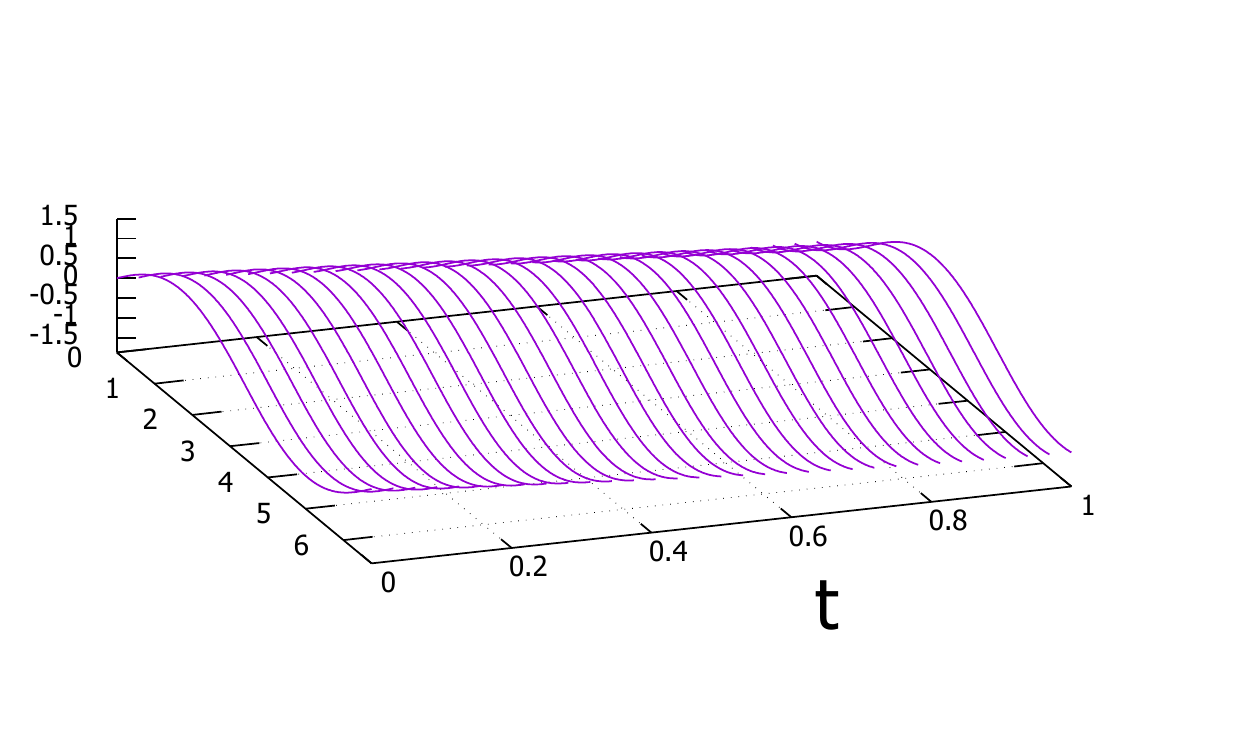}
  \includegraphics[width=7cm,clip]{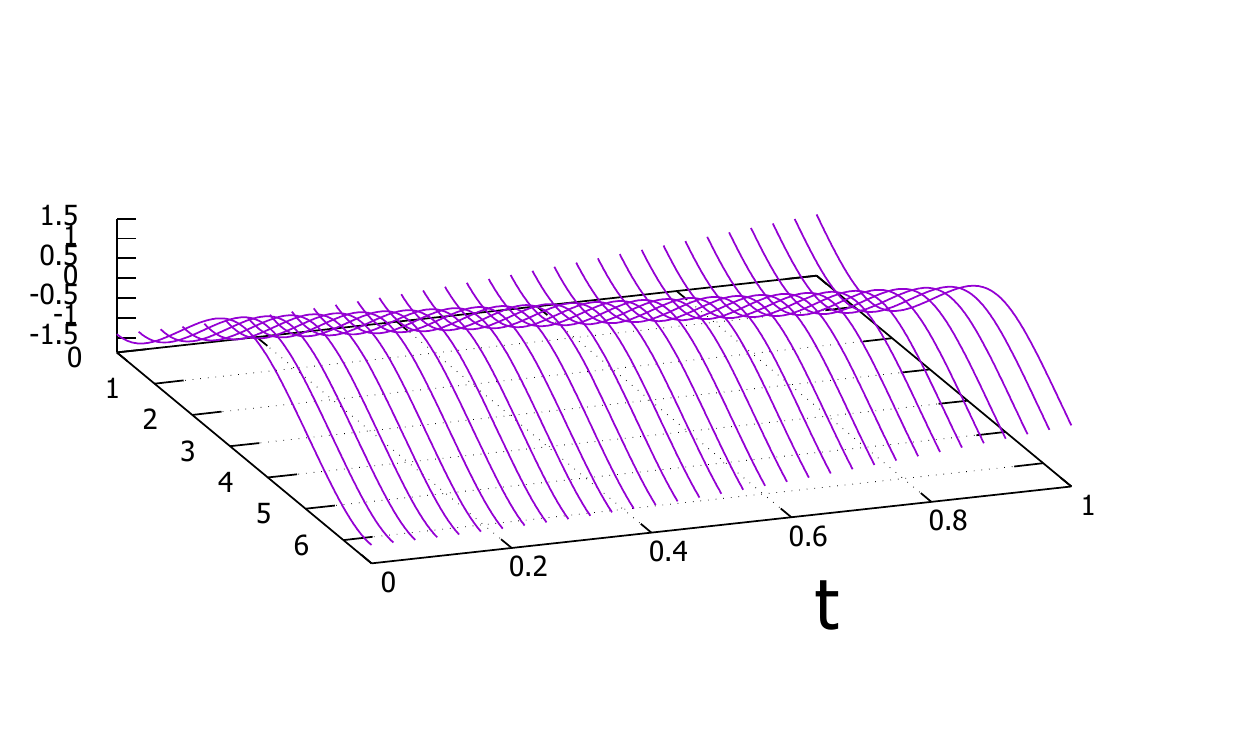}
\end{center}
\vspace{-7mm} 
\begin{minipage}{14.5cm}
  \fgcaption{Nonlinear Klein-Gordon dynamics: time evolution of $u$ and $v$.}
  \vspace{1em}
\end{minipage}
\begin{center}
  \includegraphics[width=7cm,clip]{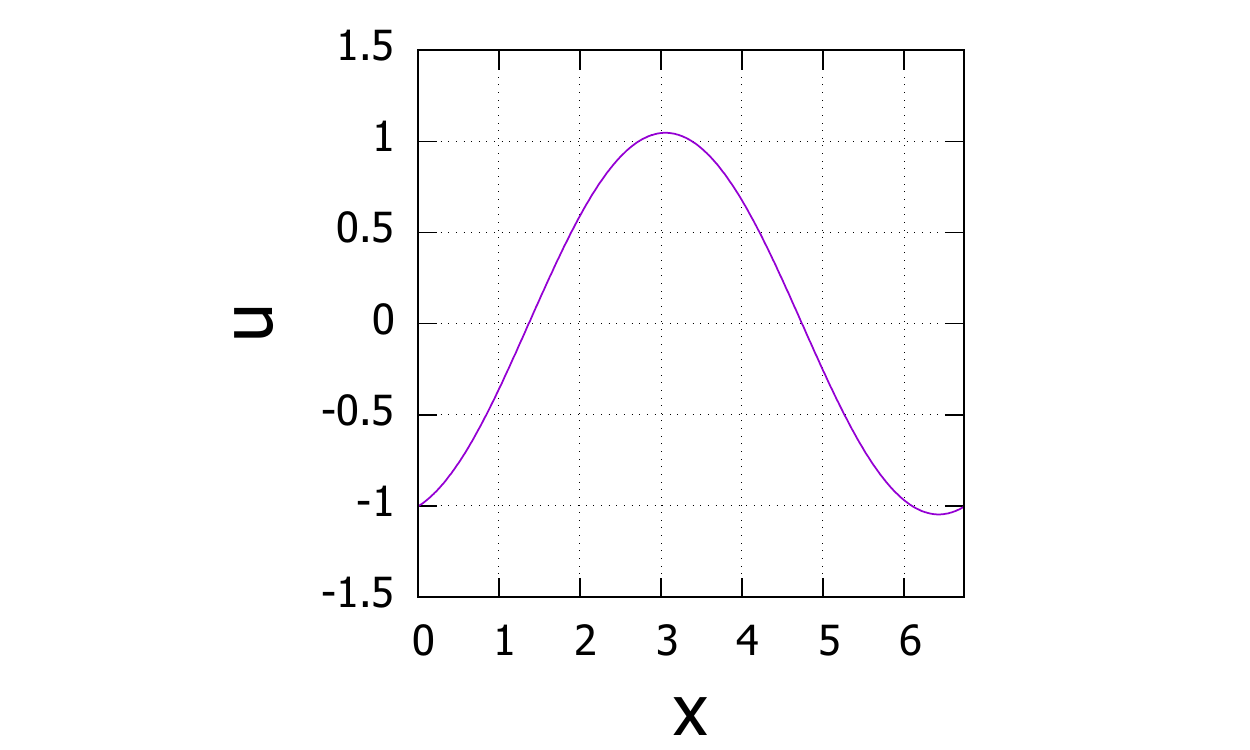}
  \includegraphics[width=7cm,clip]{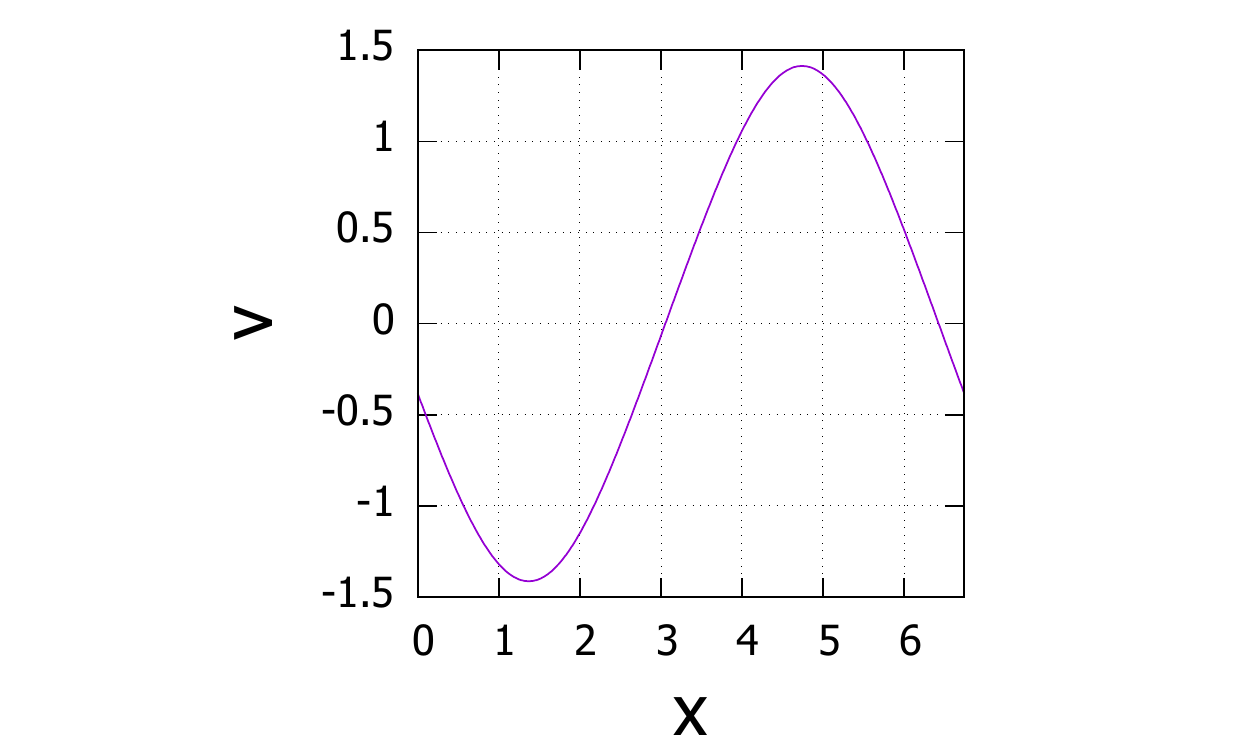}
\end{center}
\vspace{-2mm} 
\begin{minipage}{14cm}
  \fgcaption{Nonlinear Klein-Gordon equation: $u$ and $v$ at $t = 1$.}
  \label{pic:pic12a}
\end{minipage}
\end{figure*}

\subsection{Nonlinear case}
For initial and boundary values problem (2), we assume
\[ \begin{array}{ll}
F(u) = \sin u, \quad 
\alpha=-1, \quad \beta=1,  \quad
\Omega =[0,L].
\end{array}  \]
The master equation is known as the Sine-Gordon equation.
Using Jacobi’s elliptic functions  sn, cn, dn [3], the corresponding problem is written by
\begin{equation}
\label{eq:eq20}
\begin{array}{ll}
\ \tfrac{\partial v}{\partial t} -  \tfrac{\partial ^{2} u}{\partial x^{2}} + \sin u = 0, \quad
 \tfrac{\partial u}{\partial t} = v, \vspace{3mm} \\ 
\ v(x, 0) =  - \sqrt{2} \ \tfrac{\text{cn}(x, \tfrac{1}{2}) \ \text{dn}(x, \tfrac{1}{2})}{\sqrt{1 - \tfrac{1}{4} \text{sn}(x, \tfrac{1}{2})}}, \vspace{1.5mm} \\
\ v(0, t) = v(L, t), \vspace{3mm} \\ 
\ u(x, 0) = 2 \sin^{-1} \left[ \tfrac{1}{2} \text{sn} (x, \tfrac{1}{2}) \right], \vspace{3mm} \\
\ u(0, t) = u(L, t). 
\end{array}
\end{equation}
An exact solution of this problem is written by
\begin{equation}
\label{eq:eq21}
\begin{array}{ll}
\ u(x, t) = 
\ \ \ 2 \sin^{-1} \left[ \tfrac{1}{2} \ \text{sn}(x - \sqrt{2} t, \tfrac{1}{2}) \right], \vspace{3mm} \\
\ v(x, t) =
\ \ \ -  \sqrt{2} \ \tfrac{\text{cn}(x - \sqrt{2} t, \frac{1}{2}) \ \text{dn}(x - \sqrt{2} t, \tfrac{1}{2})}{\sqrt{1 - \tfrac{1}{4} \ \text{sn}^{2}(x - \sqrt{2} t, \tfrac{1}{2})}}. 
\end{array} \end{equation}
Using the perfect elliptic integral calculation, $L$ can be written by
\begin{equation}
\label{eq:eq22}
\begin{array}{ll}
L = 4 F \left( \tfrac{\pi}{2}, \tfrac{1}{2} \right)  
 \quad = 4 \int^{\pi/2}_{0} \frac{1}{\sqrt{1 - \left( \tfrac{1}{2} \right)^{2} \sin^{2} \theta }} \ d \theta \vspace{3mm} \\
 \quad =6.743001419250385098 \cdots.
\end{array}
\end{equation}
The corresponding numerical solution is shown in Figs.~7 and 8. 
In the actual numerical calculations, we further assume
\[ \begin{array}{ll}
\theta = \frac{1}{2},  \vspace{0.75em} \\ 
J \geq 2N + 1,  \vspace{0.75em} \\
L=8.
\end{array}  \]
which are exactly the same as the linear case.
The error comparing the exact and numerical solutions are shown in Figs.~7 and 8, where $N = 2^5$ and $2^{10}$ cases with $\Delta t = 2^{-2} , 2^{-3}, \cdots, 2^{-15}$ cases are examined. 
In order to obtain the error estimates, the difference between the exact and numerical solutions is calculated at each point $x_j$. 
The definition of error function follows from the linear case.
For the error arising from the time discretization, if we apply half of $\Delta t$, it results in the quartered error. This property is common to both linear and nonlinear problems. 

For the error arising from the space discretization, result in Fig.~8 shows nonlinear aspect. Indeed, error decreases depending on $N$ at first, while it can be a constant in the next. This tendency is explained by the larger effect of truncation approximation (by $N$) of Fourier expansion only in smaller $N$ cases. Indeed, the error becomes smaller if we change the value of $N$ from $2^2$ to $2^4$.
However, if we take sufficiently large $N \ge 2^5$, the error values are almost constant. 
In conclusion, even in a nonlinear case, the error of the present scheme arises mostly from the time discretization if we take sufficiently large $N \ge 2^5$.

Let us compare $N = 2^5$ to $2^{10}$ cases. 
In this comparison, there is no significant difference in error values at least if the same value for $\Delta t$ is applied. 
Note here that the convergence was turned out to be false when $N = 2^{10}$ with larger $\Delta t$ was applied. That is, the same
statement as the linear case.

\begin{figure*}[tb]
\begin{center}
  \includegraphics[width=7cm,clip]{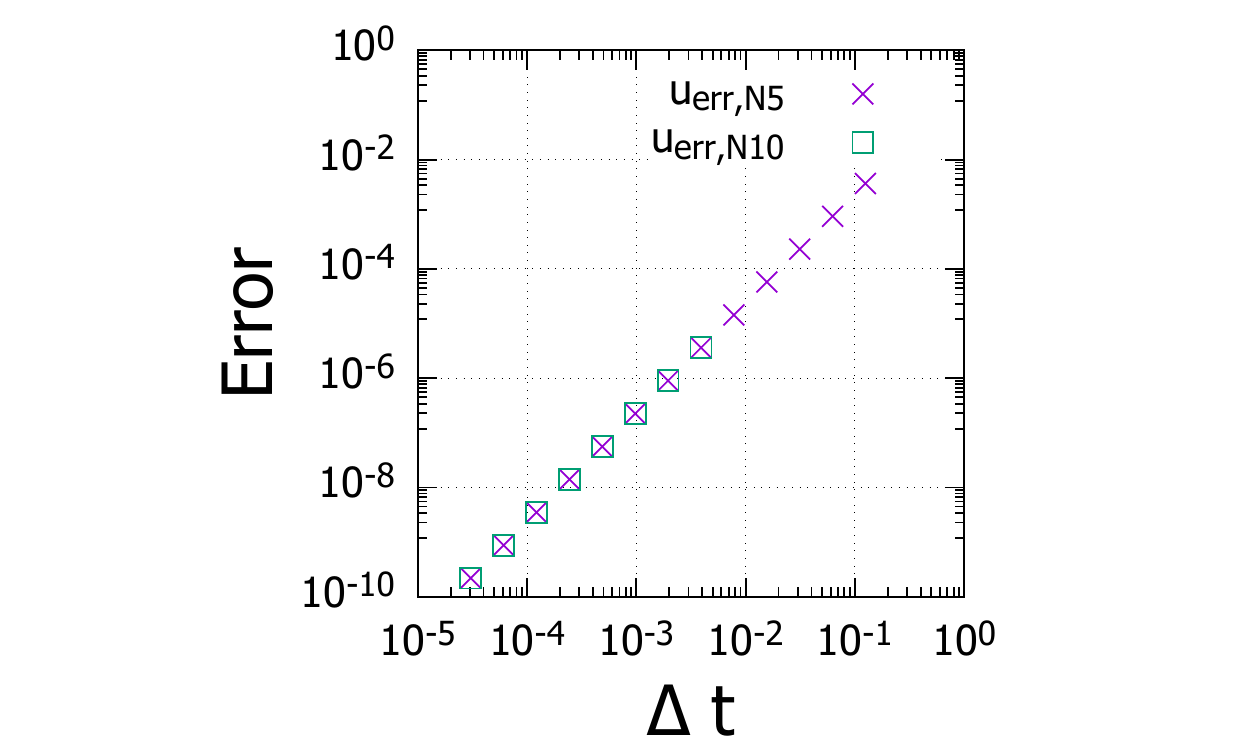}
  \includegraphics[width=7cm,clip]{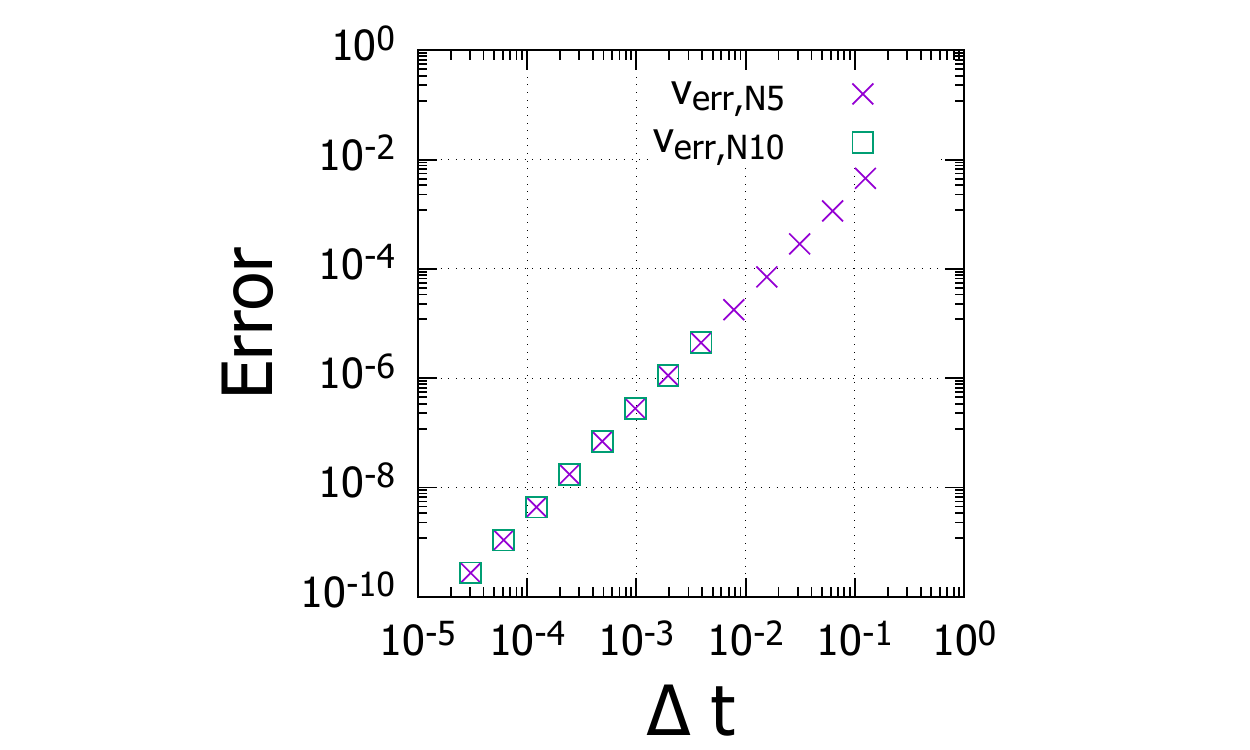}
\vspace{-2mm} 
\fgcaption{Nonlinear case: the error dependence on $\Delta t$. Two different cases $N=2^5$ (labeled by $u_{err}\_N5$ and $v_{err}\_N5$) and $N=2^{10}$  (labeled by $u_{err}\_N10$ and $v_{err}\_N10$) are examined.}
\label{pic:pic14}
\end{center}
\begin{center}
  \includegraphics[width=7cm,clip]{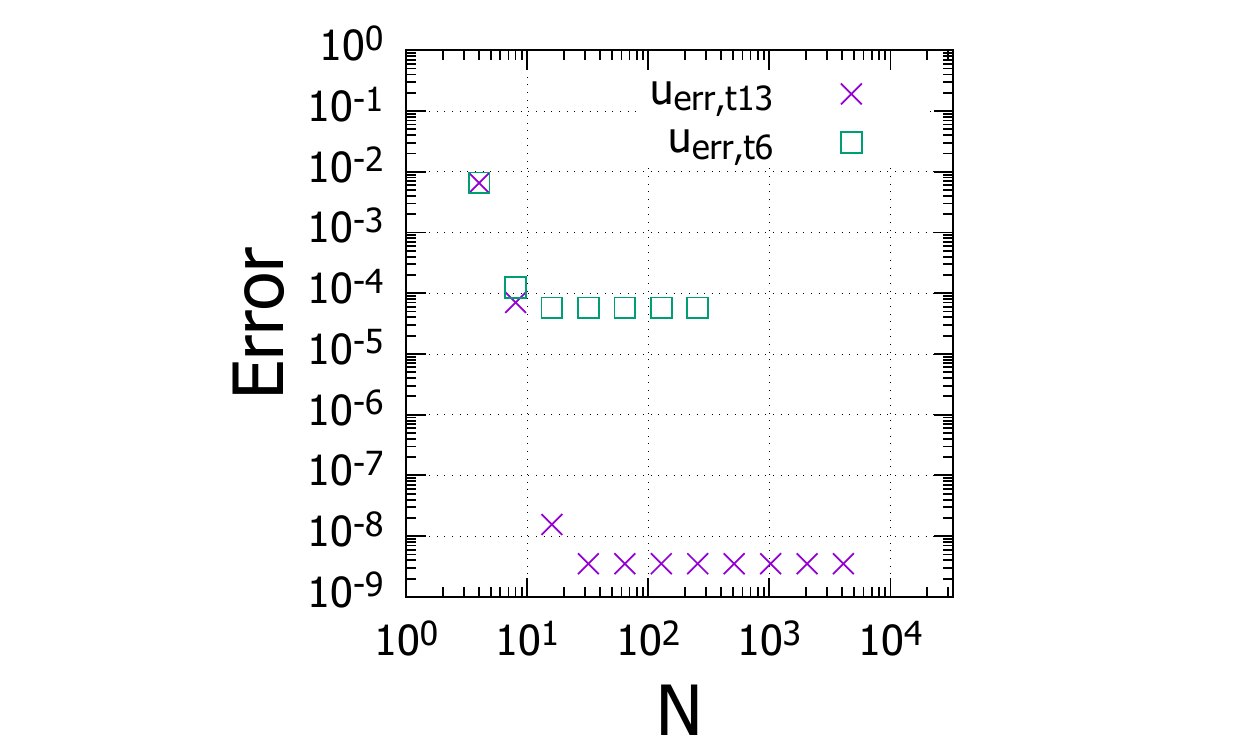}
  \includegraphics[width=7cm,clip]{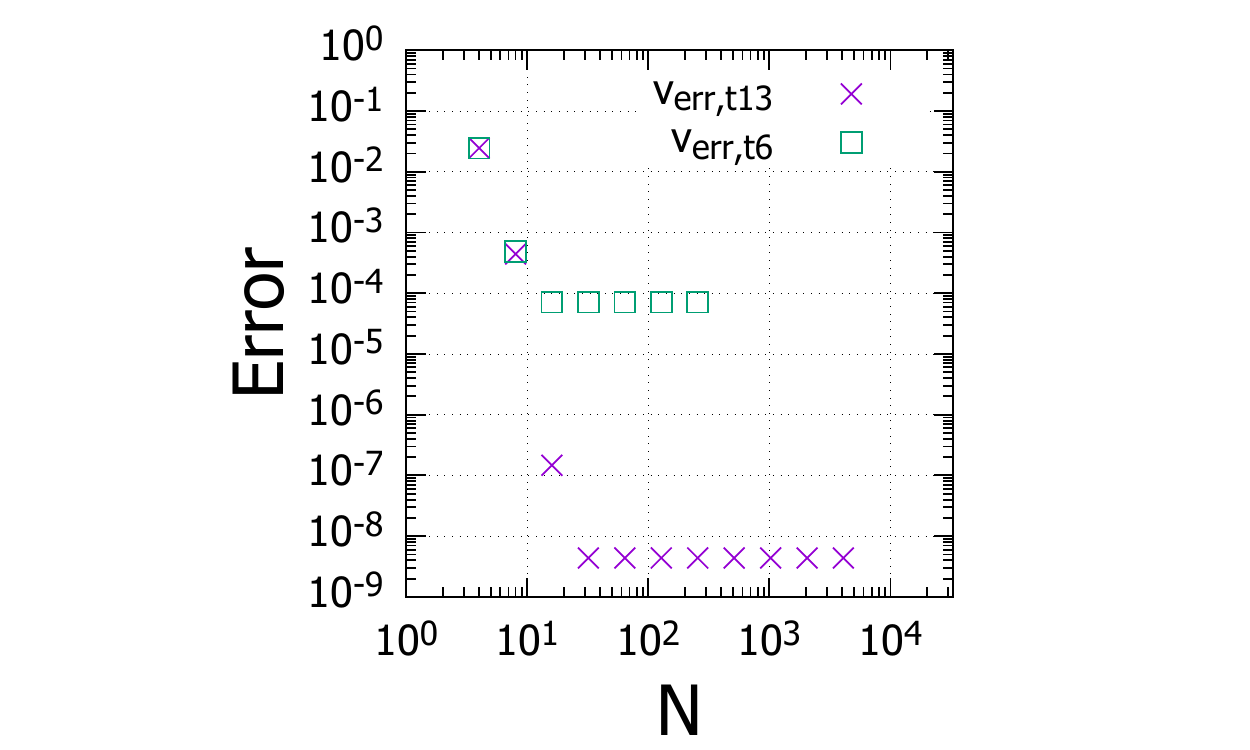}
\vspace{-2mm} 
\fgcaption{Nonlinear case: the error dependence on $\Delta t$. Two different cases $t=2^{-6}$ (labeled by $u_{err}\_t6$ and $v_{err}\_t6$) and $t=2^{-13}$  (labeled by $u_{err}\_t13$ and $v_{err}\_t13$) are examined. }
\label{pic:pic16}
\end{center}
\end{figure*}

\section{Summary}
A precise numerical scheme for nonlinear hyperbolic evolution equations is proposed; indeed, the value of Error function (corresponding to the relative error) is roughly at the order of $10^{-9}$ for a sufficiently large $N$ (Figs.~4 and 8). 
It ensures the 9-digit correctness in those benchmarks.

In terms of providing benchmark results showing the precision, the numerical solutions are compared to the exact solutions for both linear and nonlinear cases.
The relation between the numerical precision and the discretization parameters are demonstrated. 
A high precision has been confirmed to be preserved by taking sufficiently large $N \ge 2^5$, while the total calculation cost is only at the order of $N \log_2 N$. 
The error control parameters such as $N$ and $\Delta t$ are heuristically found.

Such a precise calculation would be preferably used in the wave propagation and soliton propagation in future studies. 
For some application results including the finite-dimensional representation of infinite-dimensional dynamical systems, see Ref. [4]. For the preceding works treating solitons in sub-atomic physics, see Refs.[5-7]. 
The scheme used in [5-7] is based on the finite-difference methods in which the precision is cared only at the level of obtaining convergence.

\section*{Acknowledgement}
Numerical calculations have been carried out at supercomputer at YITP, Kyoto University, and workstations at Kansai University and Tokyo Institute of Technology.
This work was partially supported by JSPS KAKENHI Grant No. 17K05440.


\end{multicols}

\end{document}